\documentclass[12pt]{article}
\usepackage{graphicx}
 \usepackage{multirow}
 \usepackage{booktabs}
 \usepackage{amsmath}
\usepackage{amsfonts}
\usepackage{amssymb, amsfonts, amsmath, nicefrac}
\usepackage{color}
\usepackage{float}
\usepackage{hyperref}
\usepackage{algorithm}
 \usepackage{algpseudocode}

\def\Gr {{\sf Gr}}

\newcommand{\cR}{{\mathfrak R}}

\newtheorem{proposition}{Proposition}[section]
\newtheorem{example}{Example}[section]
\newtheorem{remark}{Remark}[section]

\newcommand{\R}{{\cal R}}

\newcommand{\bfc}{{\bf c}}

\def\argmin{\mathop{\rm arg\,min}}

\newcommand{\U}{{\cal U}}

\newcommand{\Z}{{\cal Z}}
\newcommand{\F}{{\cal F}}

\newcommand{\be}{\begin{equation}}
\newcommand{\ee}{\end{equation}}

\newcommand{\rr}{\rightrightarrows}

\def\e{\epsilon}
\def\O{\Omega}
\def\vv{\mathfrak v}

\newcommand{\avr}{{\sf AV@R}}
\newcommand{\AVaR}{{\sf AV@R}}

\def\bbr{{\Bbb{R}}} 
\def\bbe{{\Bbb{E}}} 

\newcommand{\Ind}{{\Bbb{I}}} 

\title{\bf Risk-Averse Stochastic Optimal Control: an efficiently computable statistical upper bound}

\begin{document}

\date{ }

\maketitle

\vspace*{-2cm}

\begin{center}
\begin{tabular}{ccc}
\begin{tabular}{c}
Vincent Guigues\\
School of Applied Mathematics, FGV\\
Praia de Botafogo, Rio de Janeiro, Brazil\\
{\tt vincent.guigues@fgv.br}
\end{tabular}&
&
\begin{tabular}{c}
Alexander Shapiro\thanks{Research}\\
Georgia Institute of Technology\\
Atlanta, Georgia 30332-0205, USA,\\
{\tt ashapiro@isye.gatech.edu}\\
\end{tabular}
\end{tabular}
\end{center}

\begin{center}
\begin{tabular}{c}
Yi Cheng\\
Georgia Institute of Technology\\
Atlanta, Georgia 30332-0205, USA,\\
{\tt cheng.yi@gatech.edu}\\
\end{tabular}
\end{center}

\noindent

{\bf Abstract.}
 In this paper, we discuss an application of the
 Stochastic Dual Dynamic Programming (SDDP)   type algorithm to nested risk-averse formulations
of Stochastic Optimal
Control (SOC) problems. We
propose a construction of a statistical upper bound for the optimal value of  risk-averse SOC problems. This outlines  an approach to a solution of a  long standing problem in that area of research.
The bound holds for a large class of convex and monotone
conditional risk mappings.
Finally, we show the validity of the
statistical upper bound to solve
a real-life
stochastic hydro-thermal planning problem.
\\
\\
\noindent
{\bf Key Words:}  stochastic programming, stochastic optimal control,  SDDP, dynamic programming, risk measures,  statistical upper bounds.\\

\noindent{\bf AMS subject classifications:} 90C15, 90C90, 90C30.

\setcounter{equation}{0}
\section{Introduction}
\label{sec-intr}

Multistage stochastic optimization problems
are challenging to solve and have applications
in many areas, for instance in finance and engineering,
see for instance \cite{SDR}.
Popular methods to solve these problems
often use decomposition techniques such as Stochastic Dual Dynamic Programming (SDDP), proposed in \cite{per1991}, which is a sampling variant
of the decomposition method proposed  in \cite{bir85}.
Initially
described for risk-neutral linear problems, the
SDDP method has generated a rich literature
and many variants in the past three decades,
see, e.g., \cite{angtsoukwies19,multicutguiband,guiguescoap2013,gui2017,GSC2019, morton,cermics2018,lohnshapiro,PhilpMat,philpot8,queirozmorton}.

For risk-neutral problems and a   finite sample space, a
 stopping criterion
for SDDP is based on   estimated optimality gap determined  by deterministic lower bound
and a statistical upper bound
on the optimal value of the problem, computed   during iterations
of the method. For nested risk-averse problems, a deterministic
lower bound can be computed similar to  the risk-neutral case,
but to the best of our knowledge, no  computationally feasible
{\em statistical}
upper bound has been proposed so far for SDDP.

Of course,  in theory  the value of the constructed  approximate policy
 can be computed by  evaluating the risk at
each node of the scenario tree.
However, this computation
 rapidly becomes prohibitive
with increase of  the number   of stages and the resulting exponential growth of the  number of possible realizations of the
stochastic data process.

A deterministic upper bound on the value of
the approximate risk-averse policy was proposed in  \cite{phil13}
on the basis of inner approximations of the value functions, which is a natural extension of similar constructions for two stage programs (e.g., \cite[section 9.5]{birge-louv-book}).
Recently, two variants of Dual SDDP were introduced that also
compute a deterministic upper bound, in
\cite{cermics2018} using conjugate duality and in
\cite{GSC2019} using Lagrangian duality.
The bounds in \cite{cermics2018} and \cite{GSC2019} were developed
for risk-neutral problems, and recently extended to risk-averse problems
in \cite{dualsddpra}.
However, the computational
bulk required to compute the deterministic bounds
from \cite{phil13} and  \cite{dualsddpra} for risk-averse
problems
 increases
rapidly with increase of  the number of stages, the number of realizations
of the stochastic data per stage, and the dimension
of the state vectors.
The goal of this paper is to fill this gap
proposing an efficiently computable statistical upper bound
for SDDP applied to nested-risk averse multistage stochastic problems.
This will be possible for a large class of monotone convex risk measures that will be studied.

Our developments will be derived for Stochastic Optimal Control (SOC)
   modeling, instead of the Multistage Stochastic Programming   approach  often used in the
SDDP and related methods.    The SOC is classical with applications documented in  a large number of publications (e.g., \cite{ber78}). We would like to emphasize that many problems discussed in the    Stochastic Programming (SP) literature,  can be formulated in the SOC framework. One such example is the classical inventory model  (it is presented from both points of view, for example,  in sections 1.2.3 and 7.6.3 in \cite{SDR}). Another such example is the hydro-thermal planning problem
 discussed in section \ref{sec:numexp}.
  One  modification  in applying an  SDDP type algorithm   to SOC problems is the fact that it is not necessary anymore
to solve the dual problems to compute the required subgradients  of the cost-to-go  functions.  Of course this is a minor point since  the dual solution is often computed  by   solvers anyway.
More importantly, from the point of view of the  SDDP type algorithms,  applied to risk-averse problems,  there is an  important difference between the SOC modeling, as compared with the   SP approach.  {\em A  straightforward attempt for computation  of   statistical upper bounds in the SP framework resulted in an exponential growth of the involved  bias with increase of the number of stages, which made it practically useless }  (cf., \cite{ShapiroDing}).
On the other hand, we are going to demonstrate that in the SOC framework it is possible to construct such statistical upper bound in a computationally feasible way  for a large class of risk measures.

The outline of the paper is the following. In Section \ref{sec-spsoc},
we present the class of risk-neutral SOC problems and describe the SDDP type approach for solving
  this class of problems.
In Section \ref{preliminairies}, we present and study the risk measures
which will be used for the risk-averse SOC problem. In Section \ref{sec-raspsoc}, we present the
risk-averse SOC problem and describe the SDDP algorithm for this problem.
In Section \ref{statupperbound}, we derive our statistical upper bound.
Finally, in Section \ref{sec:numexp} we present numerical  results  where our upper bound is computed along iterations
of SDDP type algorithm  to solve a risk-averse real-life hydro-thermal planning
problem. Some additional material is given in the Appendix.

  We use the following notation. By $\xi_{[t]}:=(\xi_1,...,\xi_{t})$ we denote the history of a process
$(\xi_t)$ up to time $t$. For $a\in \bbr$, $[a]_+:=\max\{a,0\}$.
By $\Ind_A(x)$ we denote the indicator function of a set $A$, i.e., $\Ind_A(x)=0$ if $x\in A$, and  $\Ind_A(x)=+\infty$ otherwise.

\setcounter{equation}{0}
\section{Risk-neutral Stochastic Optimal Control}
\label{sec-spsoc}

Consider the   Stochastic Optimal Control  (SOC)  (discrete time, finite horizon) model (e.g., \cite{ber78}):
\begin{equation}\label{soc}
\min\limits_{\pi\in \Pi}  \bbe^\pi\left [ \sum_{t=1}^{T}
c_t(x_t,u_t,\xi_t)+c_{T+1}(x_{T+1})
\right],
\end{equation}
where $\Pi$ is the set of polices satisfying   the constraints
\begin{equation}\label{soc-b}
\Pi=\Big\{\pi=(\pi_1,\ldots,\pi_T):
u_t=\pi_t(\xi_{[t-1]}),
u_t\in \U_t, x_{t+1}=F_t(x_t,u_t,\xi_t),\;\;t=1,...,T\Big\}.
\end{equation}
Here variables  $x_t\in \bbr^{n_t}$, $t=1,...,T+1$, represent the state  of the system,   $u_t\in \bbr^{m_t}$,  $t=1,...,T$, are controls,   $\xi_t\in \bbr^{d_t}$, $t=1,...,T$, are random vectors, $c_t:\bbr^{n_t}\times\bbr^{m_t}\times\bbr^{d_t}\to \bbr$, $t=1,...,T$, are cost functions, $c_{T+1}(x_{T+1})$ is a final cost function,  $F_t:\bbr^{n_t}\times\bbr^{m_t}\times\bbr^{d_t}\to \bbr^{n_{t+1}}$ are (measurable) mappings and $\U_t$  is a (nonempty)  subset of $\bbr^{m_t}$.
Values  $x_1$  and $\xi_0$ are  deterministic  (initial conditions); it is also possible to view $x_1$ as random with a given distribution, this is not essential for the following discussion.
The optimization in \eqref{soc}   is performed over policies  $\pi\in \Pi$  determined by  decisions  $u_t$ and state variables $x_t$ considered as  functions of $\xi_{[t-1]}=(\xi_1,...,\xi_{t-1})$, $t=1,...,T$,
and satisfying the feasibility constraints \eqref{soc-b}.
  For the sake of simplicity, in order not to distract from the main message of the paper,  we assume that the control sets $\U_t$ do not depend on $x_t$.  It is possible to extend the analysis to the general case, where  the control sets are functions of the state variables,   we give a short discussion of that in section \ref{app-cont} of  the Appendix.

  With some abuse of the notation we use the same notation for $x_t$ and $u_t$, and later for $\theta_t$,  considered as functions of the random process $\xi_t$, and considered as vector   variables, e.g.,  when writing the respective dynamic programming equations. The particular meaning will be clear from the context.

 It is said that the random  process $\xi_t$ is {\em stagewise independent} if $\xi_t$ does not depend on $\xi_{[t-1]}$ for $t=1,...,T$.
\emph{}We make the following basic assumption.

\begin{itemize}
  \item [(A)]
 The random data process $\xi_1,...,\xi_T$ is stagewise independent and   its probability distribution does not depend on our decisions.
\end{itemize}
 Since it is assumed that the data process is stagewise independent, it suffices to consider policies of the form $\pi_t=u_t(x_t)$, $t=1,...,T$    (e.g, \cite{ber78}).

We can consider  problem \eqref{soc}-\eqref{soc-b} in the framework of Stochastic Programming (SP) if we view $y_t=(x_{t},u_t)$ as decision variables.
In various applications it is possible to approach the {\em same} problem using either the SOC or SP formulations.   As it was already mentioned above, for example   the classical inventory model can be treated in both frameworks (e.g., \cite[sections 1.2.3 and 7.6.3]{SDR}). Another
 such example is discussed  in section \ref{sec:numexp} below.
However, there are essential differences between the SOC and SP modeling approaches. In the SOC there is a clear separation between the state and control variables. At every stage $t$ the optimization is performed over feasible controls (also called actions)  $u_t$ and consequently the state at the next stage is determined by the state equation $x_{t+1}=F_t(x_t,u_t,\xi_t)$.
This has important  implications for the SDDP algorithm, especially in the risk averse setting.  We give a further discussion of the SOC and SP modeling approaches in Remark \ref{rem-dif} and
section \ref{app-stoch} of the Appendix.

The dynamic programming equations can be written as follows.
At the last stage, the value function
 $V_{T+1}(x_{T+1})=c_{T+1}(x_{T+1})$ and, going backward in time for $t=T,...,1$,  the value functions
\begin{equation}
\label{soc-2}
V_t(x_t)=\inf\limits_{u_t\in \U_t}
\bbe  \left [c_t(x_t,u_t,\xi_t)+
V_{t+1}\big(F_t(x_t,u_t,\xi_t) \big)\right],
\end{equation}
where the expectation is taken with respect to the (marginal) distribution of $\xi_t$,
The optimal policy is defined by the optimal controls   $\bar{u}_t(x_{t})\in\U^*_t(x_t)$, where
\begin{equation}\label{conpol-2}
\U^*_t(x_t):= \argmin_{u_t\in \U_t}
\bbe  \left [c_t(x_t,u_t,\xi_t)+
V_{t+1}\big(F_t(x_t,u_t,\xi_t) \big)\right].
\end{equation}
The optimal value of the SOC problem \eqref{soc}-\eqref{soc-b} is given by the first stage value function $V_1(x_1)$, and can be viewed  as a function of the initial conditions $x_1$.
  We make the assumptions. 
\begin{itemize}
  \item [(B)] The sets $\U^*_t(x_t)$, $t=1,...,T$,  are {\em nonempty} for every possible realization of state variables.
\end{itemize}
 
Assumption (B) holds under standard regularity conditions, e.g., if the sets $\U_t$ are compact and the objective function in the right hand side of \eqref{conpol-2} is continuous in $u_t\in \U_t$.

We consider the   convex case, by making the following assumption.
\begin{itemize}
  \item [(C)]
 For $t=1,...,T$: (i) the sets  $\U_t$   are closed convex, (ii)
 the cost functions $c_t(x_t,u_t,\xi_t)$ are convex in $(x_t,u_t)$, and
  \begin{equation}\label{affin}
 F_t(x_t,u_t,\xi_t):=A_t x_t+B_t u_t+ b_t,
 \end{equation}
with matrices
 $A_t=A_t(\xi_t)$, $B_t=B_t(\xi_t)$  and vectors $b_t=b_t(\xi_t)$ being functions of $\xi_t$.
 \end{itemize}
It follows  from Assumption (C) that the value functions $V_t(\cdot)$ are convex.
Suppose  further  that
\begin{itemize}
\item[(D)] Random vector $\xi_t$ has a finite number of realizations $\xi_{ti}$ with respective probabilities $p_{ti}$, $i=1,...,N$, $t=1,...,T$  (for the sake of simplicity assume that the cardinality $N$ is the same for every time $t$).
\end{itemize}
 
Denote $c_{ti}(x_t,u_t):=c_t(x_t,u_t,\xi_{ti})$ and $A_{ti}=A_t(\xi_{ti}),B_{ti}=B_t(\xi_{ti}),b_{ti}=b_t(\xi_{ti})$, $i=1,...,N$, the respective values of the parameters. In that case,   the    dynamic programming equations \eqref{soc-2} can be written as
\begin{equation}\label{subd-2a}
V_t(x_t)=\inf\limits_{u_t\in \U_t}
\underbrace{\sum_{i=1}^N  p_{ti} \left [c_{ti}(x_t,u_t)+
V_{t+1}\big(A_{ti} x_t+B_{ti} u_t+ b_{ti} \big)\right]}_{\displaystyle{\bbe    [c_t(x_t,u_t,\xi_t)+
V_{t+1} (A_t x_t+B_t u_t+ b_t  ) ]}}.
\end{equation}

The  subdifferentials of the value functions are  obtained from the    dynamic programming equations \eqref{subd-2a}. That is,
consider function
\[
Q_t(x_t,u_t):=\bbe  \left [c_t(x_t,u_t,\xi_t)+
V_{t+1}\big(A_t x_t+B_t u_t+ b_t \big)\right].
\]
Since $c_t(x_t,u_t,\xi_t)$ is convex in $(x_t,u_t)$ and $V_{t+1}$ is convex,
$Q_t(x_t,u_t)$  is convex.  By  \eqref{subd-2a} we have that

\begin{equation}\label{eq-q}
 V_t(x_t)=\inf_{u_t\in \U_t}Q_t(x_t,u_t)=\inf_{u_t\in \bbr^{m_t}}\left\{Q_t(x_t,u_t)+\Ind_{\U_t}(u_t)\right\}.
\end{equation}
 Consequently we have the following formula for the subdifferential of $V_t(\cdot)$ (cf., \cite[Theorem 24(a)]{roc74}):
\begin{equation}\label{roc-1}
\partial  V_t(x_t)=\big\{\gamma_t:(\gamma_t,0)\in \partial [Q_t(x_t,\bar{u}_t)+\Ind_{\U_t}(\bar{u}_t)] \big\}=
\big\{\gamma_t:(\gamma_t,0)\in \partial Q_t(x_t,\bar{u}_t)\big\},
\end{equation}
where  $\bar{u}_t$ is any point
 of $\U^*_t(x_t)$  (the indicator function can be removed in the last term of \eqref{roc-1} since the second component  of
$(\gamma_t,0)$ is 0).
  It follows that if $Q_t(\cdot,\cdot)$ is differentiable at
$(x_t,\bar{u}_t)$, then
\begin{equation}\label{roc1a}
 \nabla V_t(x_t)=\nabla Q_t(x_t,\bar{u}_t),
\end{equation}
where the gradient in the right  hand side of \eqref{roc1a}    is with respect to $x_t$.

We obtain that for any  $\bar{u}_t\in\U_t^*(x_t)$,  if functions
$c_{ti}(\cdot,\cdot)$,  $i=1,...,N$, are differentiable and
$V_{t+1}(\cdot)$ is differentiable at
$A_{ti} x_t+B_{ti} \bar{u}_t+ b_{ti}$, $i=1,...,N$, then
\begin{equation}\label{subd-3a}
\nabla V_t(x_t)=
\sum_{i=1}^N  p_{ti}  \left [\nabla  c_{ti}(x_t,\bar{u}_t)+
A^\top_{ti}  \nabla V_{t+1}\big(A_{ti} x_t+B_{ti} \bar{u}_t+ b_{ti} \big)\right].
\end{equation}
Note that a real valued convex function is differentiable almost everywhere (e.g., \cite[Theorem 25.5]{roc1970}).

Now suppose that  value functions
$V_\tau(\cdot)$
are approximated by
(lower bounding)
piecewise
affine functions
\begin{equation}\label{affine}
\underline{V}_\tau(x_\tau)=\max_{j=1,...,M} \ell_{\tau j}(x_\tau),
\end{equation}
where $\ell_{\tau j}(x_t) =a_{\tau j}^\top x_t +h_{\tau j}$, $j=1,...,M$. We need to compute a subgradient of $\underline{V}_\tau (\cdot)$  for $\tau=t+1$ when computing a subgradient of $\underline{V}_t (\cdot)$
using equation   \eqref{subd-3a}. 
A subgradient of $\underline{V}_\tau (\cdot)$ at  a point $x_\tau$ is given by $\nabla \ell_{\tau \nu}(x_t)=
  a_{\tau  \nu}$, where $\nu\in \{1,...,M\}$  is such that   $\underline{V}_\tau (x_\tau)=\ell_{\tau\nu}(x_\tau)$, i.e., $\nu$ is the index where the maximum in the right hand side of \eqref{affine} is attained and hence
  $\ell_{\tau\nu}(\cdot)$ is a supporting plane of $\underline{V}_\tau(\cdot)$ at $x_\tau$.

 This suggests a  way for computing  a subgradient of a current approximation of the value functions in a cutting planes type algorithm discussed below. There is no need to solve dual problems as in the classical SDDP method.

A cutting planes (SDDP type) algorithm for the SOC problem can be described as follows. In
the forward step at
iteration $k$ of
 the algorithm, for given convex piecewise affine lower bounding
approximations ${\underline V}_{t}^{k-1}$ of the value functions  and
 for a generated  sample path (scenario) $\hat{\xi}_1,...,\hat{\xi}_T$  of realizations of the random data process, starting with the initial value $\hat{x}_1=x_1$,   compute a minimizer in  the right hand side of \eqref{subd-2a} for the current approximation of the value function, that is
 \begin{equation}\label{subd-2ak}
\hat{u}_t\in\argmin\limits_{u_t\in \U_t}
\sum_{i=1}^N  p_{ti} \left [c_{ti}(x_t,u_t)+
{\underline V}_{t+1}^{k-1} \big(A_{ti} x_t+B_{ti} u_t+ b_{ti} \big)\right],
\end{equation}
for $x_t=\hat{x}_t$, and set  $ \hat{x}_{t+1}=F_t(\hat{x}_t,\hat{u}_t,\hat{\xi}_t)$.
If the set $\U_t$ is polyhedral and the cost functions $c_{ti}(x_t,u_t)$ are piecewise  affine  functions of $u_t$, this minimization problem  can be written as a linear programming problem, and hence has an optimal solution unless it is  unbounded from below. 
 In the next  backward step of
the algorithm,
the cutting planes approximation of the value functions are updated going backwards in time by adding the cuts at
the computed trial points $\hat{x}_t$.  These cuts
are computed using
subgradients (at the trial points) of the current approximations of the value functions.

\section{Preliminaries on risk measures}
\label{preliminairies}

Let $(\O,\F,P)$
be a probability space and let
$\Z$ be a linear  space of
$\F$-measurable functions (random variables) $Z:\O\to\bbr$.
A risk measure is a function $\R:\Z\to \bbr$
 which assigns to a
random variable $Z$ a real number representing
its risk. Typical example of the linear space $\Z$ is the space of random variables with finite $p$-th order  moments,  denoted $L_p(\O,\F,P)$, $p\in [1,\infty)$. It is said that risk measure $\R$ is {\em convex} if it possesses the properties  of
convexity, monotonicity,  and translation equivariance. If moreover it is
  positively homogeneous, then it is said that risk measure $\R$ is   {\em coherent} (coherent risk measures were introduced in
\cite{ADEH:1999}). We can refer to \cite{fol04} and \cite{SDR} for a thorough discussion of risk measures.

In this paper we consider a class of convex risk measures which can be represented in the following  parametric form:
\begin{equation}\label{coher}
 \R(Z):=\inf_{\theta \in \Theta} \bbe_P [\Psi(Z,\theta)],
\end{equation}
where $\Theta$ is a   subset of a finite dimensional vector space   and $\Psi:\bbr\times \Theta\to \bbr$ is a real valued function,  called the {\em generating function} of $\R$.  The notation $\bbe_P$   in \eqref{coher} emphasizes that the expectation is taken with respect to the probability measure (distribution) $P$ of random variable $Z$. We consider risk measures of the form \eqref{coher} for every stage. That is, for every $t=1,...,T$, we consider a   probability space  $(\O_t,\F_t,P_t)$,
and risk measure
\begin{equation}\label{coher-t}
 \R_t(Z_t):=\inf_{\theta_t \in \Theta} \bbe_{P_t} [\Psi(Z_t,\theta_t)],\;Z_t\in \Z_t,
\end{equation}
defined on  the respective  linear space of  random variables, say $\Z_t:=L_p(\O_t,\F_t,P_t)$.
For the sake of simplicity, we consider the same set $\Theta$
 and   function $\Psi$ at  every stage, this is in line with the examples below.
On the other hand, the probability distributions $P_t$ could be different for different stages.

We make the following assumptions.
\begin{itemize}
  \item [(E)]
(i) The set   $\Theta$ is  nonempty closed convex. (ii)
For every $Z_t\in \Z_t$, $t=1,...,T$,  the expectation in the  right hand side of \eqref{coher-t} is well defined and the infimum  is  finite valued.
(iii)   The function $\Psi(z,\theta)$ is convex in $(z,\theta)\in \bbr\times\Theta$. (iv)
For every $\theta \in \Theta$, the function $\Psi(\cdot,\theta)$ is   monotone nondecreasing, i.e.,
if $z_1 \leq z_2$ then
$\Psi(z_1,\theta) \leq
\Psi(z_2,\theta)$ for every $\theta\in \Theta$.
\end{itemize}

Assumption   {\rm (E)} implies that the functional $\R$, defined in \eqref{coher}, possesses the properties  of
convexity and  monotonicity.
Indeed,  it follows from assumption (E)(iii) that  $\bbe [\Psi(Z,\theta)]$ is convex in $(Z,\theta)\in \Z\times \Theta$, and
hence its minimum over convex set $\Theta$ is convex. That is, the functional $\R:\Z\to\bbr$  is convex. By Assumption  (E)(iv) the functional $\R$ is monotone, i.e., if $Z,Z'\in \Z$ are such that $Z\le Z'$ almost surely (a.s.), with respect to the measure $P$, then $\R(Z)\le \R(Z')$.


\if{
Moreover, Assumption (E)(iv) implies
that $\R$ is consistent with the stop-loss
order meaning that if $Z_1 \leq_{\mbox{icx}} Z_2$
then $\R(Z_1) \leq \R(Z_2)$
where the relation
$\leq_{\mbox{icx}}$
between random variables (representing losses) is given
by
$$
Z_1 \leq_{\mbox{icx}} Z_2
\Longleftrightarrow
\mathbb{E}[f(Z_1)]
\leq \mathbb{E}[f(Z_2)]
$$
 for any  nondecreasing convex function  $f$ such that the respective expectations are well defined
(see
\cite{sjaked02,sjaked94}).
For random variables representing
incomes, the stop-loss order
is replaced with the second order stochastic
dominance, see for instance
\cite{guiguesrom10,sjaked02,sjaked94}.
}\fi
  Recall that $Z,Z'\in \Z$ are said to be distributionally equivalent  (with respect to the reference measure $P$) if  $P(Z\le z)=P(Z'\le z)$ for all $z\in \bbr$. It is said that a functional $\R:\Z\to\bbr$ is   {\em law invariant} if $\R(Z)=\R(Z')$ for any    distributionally equivalent  $Z,Z'\in \Z$.
It follows immediately from the definition \eqref{coher-t} that $\R_t$, defined in \eqref{coher-t},
 is a function of its cdf $F_t(z)=P_t(Z_t\le z)$, and hence  is law invariant. For every $t$, consider direct product $P_1\times\cdots\times P_t$ of probability measures and the corresponding space $\Z_1\times \cdots\times \Z_t$.
  Conditional mapping   $\R_{t|\xi_{[t-1]}}:\Z_t\to \Z_{t-1}$
   is defined as a counterpart  of the law invariant functional $\R_t$, $t=1,...,T$.
  Since $\xi_0$ is  deterministic, $\R_{1|\xi_{0}}=\R$.
The associated  nested functional
is defined in the composite form
\begin{equation}\label{nested}
 \cR (\cdot):=\R_{1|\xi_0} \Big( \R_{2|\xi_{[1]}}\big (
 \cdots
\R_{T|\xi_{[T-1]}}(\cdot)\big )  \Big).
\end{equation}
 We refer to \cite[section 7.6]{SDR} for a detailed discussion of constructions of such conditional mappings and  nested functionals. Note that in this framework the process $\xi_1,...,\xi_T$,
 viewed as a random process with respect to the reference probability distributions,
 is {\em stagewise independent}  with $P_t$ being the marginal distribution of $\xi_t$.

There is a large class of   risk measures which can be represented in the parametric  form \eqref{coher}.

\begin{example}
\label{ex-avr}
{\rm
The Average Value-at-Risk measure
\begin{equation}
\label{avrmes}
\avr_{\alpha}(Z)=\inf_{\theta \in \bbr}
\bbe \left[ \theta + \alpha^{-1}[Z-\theta]_+\right],\;\alpha \in (0,1),
\end{equation}
is of form \eqref{coher} with
 generating function
$\Psi(z,\theta)=\theta+ \alpha^{-1}[z-\theta]_+$,  and $\Theta=\bbr$, $\Z=L_1(\O,\F,P)$.
  In several  equivalent  forms  the
Average Value-at-Risk was introduced over the years by different  authors in different contexts    under different names, such as Expected Shortfall, Expected Tail Loss, Conditional Value-at-Risk.    In the variational form \eqref{avrmes}  it appeared in \cite{pflug2000},\cite{ury2}.
}
 $\hfill \square$
\end{example}

\begin{example}\label{convcombexavar}
{\rm
A convex combination of
the expectation and of
Average Value-at-Risk measures is given by
\[
\R(Z):=\lambda_0 \mathbb{E} [Z] +  \sum_{i=1}^k \lambda_i \avr_{\alpha_{i}}(Z),
\]
where $\lambda_i$ are positive numbers with   $\sum_{i=0}^k \lambda_i=1$, and  $\alpha_{i}\in (0,1)$. Here $\R$  is of form
\eqref{coher} with $\Theta=\bbr^k$, $\Z=L_1(\O,\F,P)$, and
  generating function
$
\Psi(z,\theta)=\lambda_0 z + \sum_{i=1}^k \lambda_i \left(\theta_i+ \alpha^{-1}_{i}[z-\theta_{i}]_+\right).
$ } $\hfill \square$
\end{example}

\begin{example}[$\phi$-divergence]
\label{ex-div}
{\rm
Another example is  risk measures constructed from $\phi$-divergence ambiguity sets (cf., \cite{bl-2015},\cite{bental},\cite[section 7.2.2]{SDR}).
Let $\phi:\bbr\to \bbr_+\cup\{+\infty\}$ be a convex lower semicontinuous function such that $\phi(1)=0$ and $\phi(x)=+\infty$ for $x<0$. By duality arguments the distributionally robust functional associated with the ambiguity set determined by the respective $\phi$-divergence constraint with level $\e>0$ can be written in the form \eqref{coher} with
\begin{equation}\label{diverg-1}
\R_\e(Z)=\inf_{\mu,\lambda>0}\left\{\lambda\e+
\mu+\lambda\bbe_P[\phi^*((Z-\mu)/\lambda)]\right\},
\end{equation}
$\theta=(\mu,\lambda)$, $\lambda>0$, and
  generating function
$\Psi(z,\theta)=\lambda\e+
\mu+\lambda \phi^*((Z-\mu)/\lambda)$, where $\phi^*$ is the  Legendre-Fenchel  conjugate of $\phi$.
In particular for the Kullback-Leibler (KL)-divergence,
$\phi(x)=x\ln x-x+1$, $x\ge 0$, and
\begin{equation}\label{diverg-1a}
\R_\e(Z)=\inf_{\mu,\lambda>0}\left\{\lambda\e-\lambda+
\mu+\lambda \bbe_P[e^{(Z-\mu)/\lambda}] \right\}.
\end{equation}
Thus it can be represented in the form \eqref{coher}
with $\Psi(z,\lambda,\mu)=\lambda\e-\lambda+
\mu+\lambda  \,e^{(z-\mu)/\lambda}$. It could be noted that
given $\lambda>0$, the minimizer over $\mu$ in \eqref{diverg-1a}  is $\mu=\lambda\ln \bbe_P[e^{Z/\lambda}]$ and hence
\begin{equation}\label{diverg-1b}
\R_\e(Z)=\inf_{\lambda>0}\left\{\lambda\e
 +\lambda \ln\bbe_P[e^{Z/\lambda}] \right\}.
\end{equation}
However, the representation \eqref{diverg-1b} is not of the form \eqref{coher}.
} $\hfill \square$
\end{example}

Risk measures in the above  examples are positively homogeneous, and hence are coherent.

\begin{example} {\rm
Let
$u:\mathbb{R} \rightarrow [-\infty,+\infty)$ be
a proper closed concave and nondecreasing
utility function with nonempty domain.
The functional
$$
\R(Z):=\inf_{\theta \in \bbr}\big\{
\theta -\bbe[u(Z+\theta)]\big\},
$$
is of form
\eqref{coher}
with $\Theta=\bbr$ and
 generating function
$\Psi(z,\theta)=
\theta -u(z+\theta)$.
This risk measure is convex, but is not necessarily positively homogeneous. It can be viewed as  the opposite of the
OCE (Optimized Certainty Equivalent  (see \cite{bentalteboulle}).
} $\hfill \square$
\end{example}

Extended polyhedral risk measures,
introduced in \cite{guiguesrom10},  are also
of form \eqref{coher}.

\section{Risk-averse Stochastic Optimal Control}
\label{sec-raspsoc}

\subsection{Risk-averse Setting}
\label{sec-raset}
Consider  the risk averse setting in the nested form. That is, the expectation operator in the risk neutral formulation \eqref{soc} - \eqref{soc-b} is replaced by the nested risk measure $\cR$,
under the assumption that the data process is stagewise independent with respect to the reference distributions.
  Definition of $\cR$ is given in equation \eqref{nested}, and  briefly discussed in the text  above that equation.

Suppose further that the state equations are affine of the form \eqref{affin}.
This leads to the following risk averse problem (recall that $\R_{1|\xi_{0}}=\R$)   in the nested form:
\begin{equation}\label{riskav-1}
  \min_{\pi\in \Pi} \R_{1|\xi_{0}} \Big(
\bfc_1 + \R_{2|\xi_{[1]}}\big (
\bfc_2 +\cdots+
\R_{T|\xi_{[T-1]}} (\bfc_T)\big)+ \bfc_{T+1} \Big),
\end{equation}
where we use  notation $\bfc_t:=c_t(x_t,u_t,\xi_t)$, $t=1,...,T$, and
$\bfc_{T+1}:=c_{T+1}(x_{T+1})$.  The optimization (minimization) in \eqref{riskav-1} is over policies satisfying constraints \eqref{soc-b} with 
$F_t(x_t,u_t,\xi_t)$  being of the form \eqref{affin}.
The constraints  \eqref{soc-b}  should be satisfied with probability one with respect to the reference measures.   In  fact since the number of scenarios is assumed to be finite, the constraints should be satisfied for all scenarios.  Note that  as in the risk neutral case, it suffices to consider policies of the form $\pi_t=u_t(x_t)$, and that  states $x_t$ and controls $u_t$ of the considered policies are functions of $\xi_{[t-1]}$. The assumption which guarantees this is Assumption (A).

The risk averse counterpart of dynamic equations \eqref{subd-2a}
 can be written as  $V_{T+1}(x_{T+1})=c_{T+1}(x_{T+1})$ and for $t=T,...,1$,
 \begin{eqnarray}
 \label{risk-00}
 V_t(x_t)&=&\inf\limits_{u_t\in \U_t}
\R_t \big( c_{t}(x_t,u_t,\xi_t)+ V_{t+1} (A_{t} x_t+B_{t} u_t+ b_{t})\big)\\
 &=&
 \inf\limits_{u_t\in \U_t,\,\theta_t\in \Theta}  \bbe_{P_t} \left[ \Psi\big (c_{t}(x_t,u_t,\xi_t)+ V_{t+1} (A_{t} x_t+B_{t} u_t+ b_{t}),\theta_t  \big )\right],
   \label{risk-00a}
 \end{eqnarray}
where formulation  \eqref{risk-00a} is obtained by applying  definition \eqref{coher-t} of $\R_t$
with
 generating function $\Psi$.
Note that it is possible to write dynamic equations
 {\eqref{risk-00}} in terms of the (static) risk measures $\R_t$  because of the basic assumption of stagewise independence of the process $\xi_t$ (with respect to the reference measures) (e.g., \cite[section 6.5.4, Remark 39]{SDR}). The respective optimal policy $\pi_t=\bar{u}_t(x_t)$  is defined by the optimal controls
 \begin{equation}\label{risk-10}
\bar{u}_t(x_t)\in \argmin_{u_t\in \U_t}
\R_t \big( c_{t}(x_t,u_t,\xi_t)+ V_{t+1} (A_{t} x_t+B_{t} u_t+ b_{t})\big).
\end{equation}
  As in the risk neutral setting,  we assume that the set of minimizers in the right hand side of \eqref{risk-10} is {\em nonempty} for all possible realizations of state variables  (Assumption (B)).

The developments of
Section \ref{sec-spsoc}
can be   adapted
to this risk-averse
framework. Under the convexity assumption  (C), the value functions $V_t(\cdot)$ are convex in the risk averse setting as well. There are explicit
formulas
 how to compute a subgradient of the functional $\R:\Z\to \bbr$  for various examples of  risk measures  (cf., \cite[section 6.3.2]{SDR}).

Recall definition \eqref{coher-t} of risk measure $\R_t$. For $x_t$ and the  optimal control $\bar{u}_t=\bar{u}_t(x_t)$, determined by \eqref{risk-10}, consider a minimizer
 \begin{equation}\label{minimizer}
 \bar{\theta}_t\in\argmin_{\theta_t\in \Theta} \bbe_{P_t} \left[ \Psi\big (c_{t}(x_t,\bar{u}_t,\xi_t)+ V_{t+1} (A_{t} x_t+B_{t} \bar{u}_t+ b_{t}),\theta_t  \big )\right].
 \end{equation}
  Note that $\bar{\theta}_t$ can be computed in two equivalent ways. One way is to solve the minimization problem \eqref{risk-00a} jointly in $u_t$ and $\theta_t$. The other approach is to use \eqref{minimizer} using computed optimal controls $\bar{u}_t$. In that case $\bar{\theta}_t$ is a function of $\bar{u}_t$
which in turn is a function of $\xi_{[t-1]}$.
In both cases $\bar{\theta}_t$ can be viewed as a function of $\xi_{[t-1]}$. In the following developments we use the second approach since it is relatively easy to compute $\bar{\theta}_t$ using formula \eqref{minimizer}.

Then, similar to \eqref{subd-3a} and using the
Chain rule, a subgradient $\nabla V_t(x_t)$ of the value function $V_t$  at $x_t$ can be computed as
 \begin{equation}\label{subdif-3ra}
\nabla V_t(x_t)=
\bbe_{P_t}  \left [
  \Psi'(y_t
,\bar \theta_t)
\Big(
\nabla   c_t(x_t,\bar{u}_t,\xi_t)+
A^\top_t  \nabla V_{t+1}\big(A_t x_t+B_t \bar{u}_t+ b_t \big)\Big)\right],
\end{equation}
where
$\Psi'(y_t,\bar \theta_t)$
is a subgradient  of
$\Psi(\cdot,\bar \theta_t)$
at $y_t$, $\nabla   c_t(x_t,\bar{u}_t,\xi_t)$ is a
subgradient of $c_t(\cdot,\bar{u}_t,\xi_t)$ at $x_t$,
$\nabla   V_{t+1}(
A_t x_t + B_t \bar u_t + b_t)$
is a
subgradient of
$V_{t+1}$ at
$A_t x_t + B_t \bar u_t + b_t$,
and
$y_t:=c_{t}(x_t,\bar u_t,\xi_t)+ V_{t+1} (A_{t} x_t + B_{t} \bar u_t+ b_{t})$.
(If
$\Psi(\cdot,\bar \theta_t)$ is differentiable at $y_t$, then
$\Psi'(y_t,\bar \theta_t)$ is given by the derivative of $\Psi(\cdot,\bar \theta_t)$ at $y_t$.)

 As a special case, consider  Example  \ref{ex-avr} of the Average Value-at-Risk measure. In that case the minimizer $\bar{\theta}$ in the right hand side of \eqref{avrmes} is given by the $(1-\alpha)$-quantile of the considered distribution. That is,  suppose that the reference distribution $P_t$ has a finite number of $N$ realizations  with equal
probabilities $1/N$. Then  $\bar{\theta}_t$ can be computed by arranging values $c_{ti}(x_t,\bar{u}_t)+ V_{t+1} (A_{ti} x_t+B_{ti} \bar{u}_t+ b_{ti})$,  $i=1,\ldots,N$,  in the increasing order and taking the respective empirical $(1-\alpha)$-quantile. Consequently, the required subgradient of the current   lower  approximation of the value function can be computed in a straightforward way (cf., \cite{Sha2012a}).

\subsection{Statistical upper bounds on the value of the policy}\label{statupperbound}

 In this section, we discuss the construction of
a statistical
upper bound on the optimal  value of the risk averse problem.  As before, all probabilistic statements and expectations are taken with respect to the {\em reference distributions}.
Let $\underline{V}_{t}(x_t)$, $t=1,...,T$,  be   current approximations of the value functions.
This defines the corresponding (approximate) policy
 $(\hat x_t, \hat u_t)$ with
  \begin{equation}\label{aprpol-1}
\hat{u}_t \in \argmin_{u_t\in \U_t}
\R_t \big( c_{t}(\hat x_t,u_t,\xi_t)+
 \underline{V}_{t+1} (A_{t} \hat x_t+B_{t} u_t+ b_{t})\big).
\end{equation}
Observe that by the construction, $V_t(\cdot)\ge \underline{V}_{t}(\cdot)$ for $t=1,...,T$, and hence
value $\underline{V}_{1}(x_1)$ gives a lower bound for the optimal value of the considered  problem.

 For a given realization (scenario) $\xi_1,...,\xi_T$ of the data process,  $\hat x_t$ and  $\hat u_t$ are computed in the forward step of the SDDP algorithm, and can be viewed as functions
 $\hat x_t=\hat x_t(\xi_{[t-1]})$ and
 $\hat u_t= \hat u_t (\xi_{[t-1]})$. When   each reference probability distribution has a finite support (of $N$ points), i.e., for  the discretized  version of the problem, these  values are computable.

Now let $\hat{\theta}_t\in \Theta$ be a specified function of the data process, $\hat{\theta}_t=\hat{\theta}_t(\xi_{[t-1]})$, $t=1,...,T$.
 Note that   $\hat{\theta}_t$ is non-anticipative in the sense that it does not depend on unobserved values $\xi_t,...,\xi_T$ at time $t$.
 Denote $\hat c_t:=c_{t}(\hat x_{t},\hat u_t,\xi_t)$, $t=1,...,T$, and
 $\hat c_{T+1}:=c_{T+1}(\hat x_{T+1})$.
 Consider the following sequence of random variables (functions of the data process)  defined iteratively going backward in time:   $\vv_{T+1}:=\hat c_{T+1}$
and
\begin{equation}\label{recur}
 \vv_{t}:=\Psi(\hat c_t + \vv_{t+1}, \hat \theta_t),\;t=T,\ldots,1.
\end{equation}
 Of course, values $\vv_{t}$ depend on a choice of parameters $\hat \theta_t$.
 We will discuss an appropriate choice of $\hat{\theta}_t$ later.
 Our statistical upper bound
on the value of a risk-averse approximate policy is given in the following proposition.

\begin{proposition}
\label{pr-value}
Consider the  risk-averse   problem \eqref{riskav-1}
 Let $\vv_t$ be the sequence of random variables  (defined iteratively by \eqref{recur}) associated with   current approximations of the value functions. Then  for $t=1,...,T$,
\begin{equation}\label{uppbound}
 \R_{t|\xi_{[t-1]}}\big (
\hat c_{t}+\ldots
  +
\R_{T|\xi_{[T-1]}}( \hat c_T + \hat c_{T+1}) \big)\le \bbe_{|\xi_{[t-1]}}[\vv_t],\;\;w.p.1.
\end{equation}
In particular,  $\bbe[\vv_1]$ is greater than or equal to  the   value of the policy defined  by the considered approximate value functions,  and
is an upper bound on the optimal  value of the
 risk averse problem.
\end{proposition}

  {\bf{Proof.}}
For $t=T$,
using the definition of
$\hat u_T$ and since $\hat\theta_T\in \Theta$,
  we get
$$
\begin{array}{lcl}
\R_{T|\xi_{[T-1]}} (\hat c_T + \hat c_{T+1})
&=& \inf\limits_{u_T\in \U_T}
\R_T \left(c_{T}(\hat{x}_T,u_T,\xi_T)+
\hat{V}_{T+1} (A_{T}
\hat{x}_T+B_{T} u_T +b_{T})\right)\\
&\le&
\bbe_{|\xi_{[T-1]}} \Big[ \Psi\Big(c_{T}(\hat{x}_T,\hat u_T,\xi_T)
+c_{T+1}(A_T \hat x_T + B_T \hat u_T+b_T),\hat\theta_T\Big) \Big]\\
&=&\bbe_{|\xi_{[T-1]}}  [ \vv_T ].
\end{array}
$$
We now use induction in $t$ going backward in time.
For $t-1$  we have
\begin{equation}\label{proofinductionub1}
\begin{array}{l}
\R_{t-1|\xi_{[t-2]}} \Big(
\hat c_{t-1} + \R_{t|\xi_{[t-1]}}\big(
\hat c_{t}+\ldots
  +
\R_{T|\xi_{[T-1]}}(\hat c_T + c_{T+1}(\hat x_{T+1}))\big)\Big) \\
\leq
\R_{t-1|\xi_{[t-2]}} \big (
\hat c_{t-1}+ \bbe_{|\xi_{[t-1]}}[\vv_t] \big) \;\;\;\mbox{(monotonicity and induction step)} \\
\leq
\bbe_{|\xi_{[t-2]}}\big[
\Psi\big(\hat c_{t-1}+ \mathbb{E}_{|\xi_{[t-1]}}[\vv_t],\hat \theta_{t-1}\big)\big]\;\;  \mbox{ (because ${\hat \theta}_{t-1} \in \Theta$)}\\
=\bbe_{|\xi_{[t-2]}}\big[
\Psi\big(\bbe_{|\xi_{[t-1]}}[\hat c_{t-1}+\vv_t],\hat \theta_{t-1}\big)\big]\;\;  \mbox{(since $\hat c_{t-1}$
is a function of $\xi_{[t-1]}$)}
\\
\le  \bbe_{|\xi_{[t-2]}}
\bbe_{|\xi_{[t-1]}}\big[\Psi\big(\hat c_{t-1}+\vv_t,\hat \theta_{t-1}\big)\big]   \mbox{(by  Jensen's inequality})\\
=\bbe_{|\xi_{[t-2]}}
 \big[\Psi\big(\hat c_{t-1}+\vv_t,\hat \theta_{t-1}\big)\big]\\
 =\bbe_{|\xi_{[t-2]}}[\vv_{t-1}].
\end{array}
\end{equation}
This completes the induction step.
$\hfill \square$
\\

Therefore, for a sample path (scenario) of the data process, an unbiased point estimate of an upper bound on the corresponding policy  value can be computed recursively
 starting with $\vv_{T+1}=c_{T+1}(\hat{x}_{T+1})$ and
  going backward in time using the iteration procedure \eqref{recur}.
  Finally $\vv_1$ gives   a  point estimate of an upper bound
  on the corresponding value of the policy.
Therefore  by generating a sample of scenarios, of the random data process,  and averaging the corresponding point estimates it is possible to construct   the respective statistical upper bound for the optimal value of the risk averse problem.

The quality of such statistical bound depends on the choice of the parameter  function $\hat\theta_t$. It is natural to use the corresponding  minimizer of the form \eqref{minimizer}. That is, to take
\begin{equation}\label{minimizer-2}
 \hat{\theta}_t\in\argmin_{\theta_t\in \Theta} \bbe  \left[ \Psi\big (c_{t}(\hat x_t,\hat{u}_t,\xi_t)+  \underline{V}_{t+1} (A_{t} \hat x_t+B_{t} \hat{u}_t+ b_{t}),\theta_t  \big )\right].
 \end{equation}
 The so defined  $ \hat{\theta}_t$ is a function of $\hat x_t$ and $\hat{u}_t$, which in turn are functions of $\xi_{[t-1]}$.
For example, as it was pointed at the end of Section \ref{sec-raset}, in case of the Average Value-at-Risk measure such  $\hat{\theta}_t$ can be easily computed by using  the respective quantile. Note that even for $ \hat{\theta}_t$ of the form \eqref{minimizer-2} the inequality \eqref{uppbound} can be strict. This is because Jensen's inequality was used in derivations \eqref{proofinductionub1}.  Nevertheless, this approach performed well in the numerical experiments discussed in the next section.

\begin{remark}
\label{rem-dif}
{\rm
We would like to point to the important difference between the corresponding  SOC and SP approaches to construction of the statistical upper bound for the risk averse problems.
Computation of the  parameter
$ \hat{\theta}_t $ in \eqref{minimizer-2} is   based on  the distribution of random vector $\xi_t$. When  $\xi_t$ has a finite  number of realizations $\xi_{ti}$, $i=1,...,N$,  the parameter  $ \hat{\theta}_t $ is a function of   all corresponding costs $\hat{c}_{ti}$ and   {\em all} values   $A_{ti},B_{ti},b_{ti}$, $i=1,...,N$, of random parameters at stage $t$. This makes  $ \hat{\theta}_t $,     in a sense,  to be  a ``consistent" estimate of $\bar{\theta}_t$ defined  in \eqref{minimizer}.  On the other hand, in the SP setting it was not possible to construct  a computationally feasible    consistent estimate of the respective parameter of the risk measure.  As a result a straightforward attempt for computation  of such statistical upper bound in the SP framework resulted in an exponential growth of the involved  bias with increase of the number of stages, which made it practically useless (cf., \cite{ShapiroDing}). $\hfill \square$
}
\end{remark}

  We close this section by presenting  Algorithm \ref{alg:SDDP-SOC} for computing the statistical upper bound for a $T$-stage SOC problem.
\begin{algorithm}
	\caption{SDDP-type Algorithm for SOC Problem}
	\begin{algorithmic}[1]
		\State Inputs: stage-wise independent samples ${\xi}_t :=\{{\xi}_{tj}\}_{1\le j \le N_t}, t=1,\cdots,T,$ initializations of $V_t(\cdot): \underline{V}^0_t(\cdot), t =1,\cdots, T,$ initial point $\hat{x}_1$
		\For {$k=1,2,\ldots,K$}
		\State
		$\underline{V}^{k-1}_{T+1}(\cdot)=V_{T+1}$
		\For {$t = 1,\cdots, T$}\Comment{Forward Step}
		\State $\hat{u}_t = \argmin\limits_{u_t\in \mathcal{U}_t} \mathcal{R}_t\left( c_t(\hat{x}_t,u_t,{\xi}_t) + \underline{V}_{t+1}^{k-1}(A_t\hat{x}_t + B_t u_t + b_t)\right)$
		\State Draw a sample $(\hat{A}_t, \hat{B}_t, \hat{b}_t)$ from $\{{\xi}_t\}$
		\State $\hat{x}_{t+1} = \hat{A}_t\hat{x}_t + \hat{B}_t\hat{u}_t + \hat{b}_t$
		\EndFor
		\For {$t = T,\cdots, 1$}\Comment{Backward Step}
		\State {$\hat{\theta}_t = \argmin\limits_{\theta_t\in \Theta} \frac{1}{N}\sum\limits_{j=1}^{N} \Psi\left(c_t(\hat{x}_t, \hat{u}_t, \xi_{tj}) + \underline{V}^{k-1}_{t+1}(A_{tj}\hat{x}_t + B_{tj}\hat{u}_t + b_{tj}),\theta_t \right)$,}  \label{alg-compute-theta-backward}
		\State $v_t = \frac{1}{N}\sum\limits_{j=1}^{N} \Psi\left( c_t(\hat{x}_t, \hat{u}_t, \xi_{tj}) + \underline{V}^{k-1}_{t+1}(A_{tj}\hat{x}_t + B_{tj}\hat{u}_t + b_{tj}),\hat{\theta}_t\right),$
				\State $y_{t j}:=c_{t}(\hat x_t,\hat u_t,\xi_{t j})+ \underline{V}_{t+1}^{k-1}(A_{t j} \hat x_t + B_{t j} \hat u_t+ b_{t j})$,
		\State $g_t = \frac{1}{N}\sum\limits_{j=1}^{N} \Psi'(y_{t j},\hat{\theta_t})\left(\nabla c_t(\hat{x}_t, \hat{u}_t,{\xi}_{tj}) + {A}^{\top}_{tj}\nabla \underline{V}_{t+1}^{k-1}({A}_{tj}\hat{x}_t + {B}_{tj}\hat{u}_t + {b}_{tj})\right),$
		\State $\underline{V}^k_t(x_t)=\max(\underline{V}^{k-1}_t(x_t),g_t^T (x_t - \hat{x}_t) + v_t)$,
		\EndFor
		\State Lower bound: $L_k = \underline{V}^k_1(\hat{x}_1)$
		\State Generate $S$ sample paths $\xi_s^k=\{\xi_{ts}^k\}_{1\le t \le T}, s=1, \cdots, S$, run forward step for each sample  $\vphantom\;\;\;\;\;\;$path $\xi_{s}^k$ to obtain controls $(\hat{u}_{ts}^k)_{1\le t \le T}$ and states $(\hat{x}_{ts}^k)_{1\le t \le T+1}$ \Comment{Evaluation}
		\State Set $\vv_{T+1, s}^k = c_{T+1}(\hat{x}_{T+1, s}^k), s=1,\cdots, S$
		\For{$t = T,\cdots, 1$}
		\For{$s = 1,\cdots, S$}
  \State $\hat{\theta}_{t s}^k = \argmin\limits_{\theta_t\in \Theta} \frac{1}{N}\sum\limits_{j=1}^{N} \Psi\left(c_t(\hat{x}_{t s}^k, \hat{u}_{t s}^k, \xi_{t j}) +  \underline{V}^{k}_{t+1}(A_{t j} \hat{x}_{t s}^k + B_{t j} \hat{u}_{t s}^k + b_{t j} ),\theta_t \right)$  \label{alg-compute-theta-eval}
		\State  $\vv_{ts}^k = \Psi(c_t(\hat{x}_{ts}^k,\hat{u}_{ts}^k,\xi_{ts}^k)+\vv_{t+1,s}^k,\hat{\theta}_{t s}^k)$
		\EndFor
		\EndFor
		\State $\bar{\vv}_1^k = \frac{1}{S}\sum\limits_{s=1}^S \vv_{1s}^k, \sigma_k^2 = \frac{1}{S-1}\sum\limits_{s=1}^S(\vv_{1s}^k - \bar{\vv}_1^k)^2$
		\State Statistical upper bound: $U_S^k = \bar{\vv}_1^k + z_{1-\beta}\sigma_k/\sqrt{S}$.
\EndFor
	\end{algorithmic}
	\label{alg:SDDP-SOC}
\end{algorithm}

\section{Numerical Experiments}\label{sec:numexp}

 In this section  numerical experiments are performed on the Brazilian Inter-connected Power System problem (we refer to \cite{Sha2012a} for more details on the problem description). All experiments were run using Python 3.8.5 under Ubuntu 20.04.1 LTS operating system with a 4.20 GHz Intel Core i7 processor and 32Gb RAM.    We extended the \texttt{MSPPy} solver \url{{https://github.com/lingquant/msppy}}
 \cite{ding2019} for the SDDP algorithm solving for the SOC problem.  We report numerical results of the convergence guided by the deterministic lower bound and the statistical upper bound of the risk averse stochastic optimal control problem.

The hydro-thermal planning problem is a large-scale problem with  $T=120$  planning horizon stages  and four state variables related to the energy reservoirs  in four interconnected regions. The monthly energy inflows define the  stochastic  data process in the model. For the sake of simplicity, it is assumed in the experiments below that the random inflow   process is stagewise independent. The (discretization) samples are generated from log-normal distributions (with $100$ realizations at each stage) estimated from the historical data. Previous attempts to define a statistical upper bound have shown some of the challenges of this task. For example, the numerical results in \cite{ShapiroDing} show that by formulating the problem as a risk-averse multistage stochastic program, the scale of the statistical upper bounds starts to explode
with increase of the number of stages and becomes
prohibitively large
when the number of stages $T$ is more than $10$.

We aim to demonstrate via the hydro-thermal planning problem, the effectiveness of the construction of the statistical upper bound proposed  in Section \ref{sec-raspsoc}. This  suggests first to formulate the problem as a risk-averse optimal control model, and then to solve it by a variant of the SDDP algorithm, while preserving the number of stages,  the states, and the data process in the original problem. More specifically, we construct the upper bound as explained in Section \ref{statupperbound}, detailed in Algorithm \ref{alg:SDDP-SOC}. We conduct experiments for risk measures of convex combination of expectation and \AVaR \ and KL-divergence, as described in Examples \ref{convcombexavar} and   \ref{ex-div}, respectively. We solve both problems, and compute the corresponding statistical upper bounds,   by an SDDP-type algorithm as described  in Algorithm \ref{alg:SDDP-SOC}.

\paragraph{Implementation Details.}
\begin{enumerate}
	\item Convex combination of expectation and $\AVaR$ (Example \ref{convcombexavar}): $(1-\lambda)\mathbb{E}[\cdot] + \lambda \avr_{\alpha}(\cdot)$. For this risk measure, we perform  tests with
   the  critical value of the confidence interval $z_{1-\beta}= 2$ (see line 26 of Algorithm \ref{alg:SDDP-SOC})  and  $\lambda \in \{0, 0.5, 1\}$. When $\lambda = 0$, the problem becomes  risk neutral, while $\lambda = 1$ corresponds to an extreme risk aversion.

In this setting, at each backward step and in the evaluation procedure (line \ref{alg-compute-theta-backward} and line \ref{alg-compute-theta-eval} in Algorithm \ref{alg:SDDP-SOC}), $\hat{\theta}_t$ can be computed by arranging values $c_t(\hat{x}_t,\hat{u}_t,\xi_{tj}) + \underline{V}_{t+1}(A_{tj}\hat{x}_t + B_{tj}\hat{u}_t + b_{tj}), j = 1,\cdots, N$, in the increasing order and taking the respective empirical $(1 - \beta)$-quantile. Moreover, in order to obtain a fast converging deterministic lower bound, we adopt the biased-sampling technique proposed in \cite{rasddpliushap20}.

 \item KL-divergence (Example \ref{ex-div}). For this risk measure, we conduct experiments for  $\epsilon \in \{10^{-1}, 10^{-2}, 10^{-3}, 10^{-8}, 10^{-12}\}$, which corresponds to problems with different levels of risk aversion. In particular, when $\epsilon = 10^{-12}$, the problem is essentially a risk neutral problem, up to some numerical error.

 In this case, at steps indicated by line \ref{alg-compute-theta-backward} and line \ref{alg-compute-theta-eval} in Algorithm \ref{alg:SDDP-SOC}, the following (one-dimensional) convex program:
 \begin{equation}
 	\hat{\lambda}_t = \argmin\limits_{\lambda_t > 0} \{\lambda_t\epsilon + \lambda_t\ln\mathbb{E}_{P_t}\left[ e^{Z_t/\lambda_t} \right]\},
 \end{equation}
 where $Z_t := \{c_t(\hat{x}_t,\hat{u}_t,\xi_{tj}) + \underline{V}_{t+1}(A_{tj}\hat{x}_{t} + B_{tj}\hat{u}_t + b_{tj})\}_{1\le j \le N_t}$,
was solved using
\texttt{Scipy} solver.

\end{enumerate}

\paragraph{Results.}
For risk measure $(1-\lambda)\mathbb{E}[\cdot] + \lambda \avr_{\alpha}(\cdot)$ , with  $\lambda = 0.5$,
in order to examine the trend of the statistical upper bound, we compute the upper bound for the problem   at every 10 iterations with a sample of size $S=10$, by running 10 forward passes in parallel.
Figure \ref{fig:evo} in the Appendix displays the evolution of the deterministic lower bounds and the statistical upper bounds for the hydro-thermal planning problem for 3000 iterations. We can see from the figure that the statistical upper bound oscillates significantly for the first 500 iterations and then gradually stabilizes within narrow fluctuations. Table \ref{tab:mc}
reports,  for different choices of $\lambda$, the statistical upper bounds obtained from Monte Carlo simulation using
3000 samples, along with the deterministic lower bounds and the relative gap ($\frac{\text{upper bound}\ - \text{lower bound}}{\text{lower bound}}$) at
the last
iteration 3000. From the results, it seems that the relative gap of the problem is not very sensitive to the level of risk aversion.

\begin{table}[h!]
\centering
 		\begin{tabular}{cccc}
 			\toprule
 			\multicolumn{4}{c}{$(1-\lambda)\mathbb{E}[\cdot] + \lambda \avr_{\alpha}(\cdot)$}\\
 			\midrule
 			\multicolumn{1}{l}{ $\lambda$} & \multicolumn{1}{l}{Deterministic lower bound} & \multicolumn{1}{l}{Statistical upper bound} & \multicolumn{1}{l}{Gap($\%$)} \\
 			& ($\times 10^{9}$) & ($\times 10^{9}$)\\
 			\midrule
 			0.0         & 0.345  &0.348 &0.97 \\
 			\\[-1em]
 			0.5         &1.640       &1.672     &1.93     \\
 			\\[-1em]
  			1.0        &6.669       &7.003     &5.02  \\
 			\\[-1em]
 			\bottomrule
 		\end{tabular}
 		\caption{Convergence of convex combination of expectation and \AVaR \ problem for different $\lambda$.}
 		\label{tab:mc}
\end{table}

 Table \ref{tab:mc2} reports results  for the KL-divergence problem. The statistical upper bounds are computed by Monte Carlo simulation using $3000$ samples,    the lower bound and the relative gap, are computed as well  for difference values  of $\epsilon$. All  results in the table are obtained when the problems are solved for $3000$ iterations. We observe   that when  $\epsilon$ increases, the relative gap becomes larger.

\begin{table}[h]
\centering
 		\begin{tabular}{cccc}
 			\toprule
 			\multicolumn{4}{c}{KL-divergence}\\
 			\midrule
 			\multicolumn{1}{l}{ $\epsilon$} & \multicolumn{1}{l}{Deterministic lower bound} & \multicolumn{1}{l}{Statistical upper bound} & \multicolumn{1}{l}{Gap($\%$)} \\
 			& ($\times 10^{9}$) & ($\times 10^{9}$)\\
 			\midrule
 			$10^{-1}$         &4.894    &5.959   &21.76  \\
 			\\[-1em]
 			$10^{-2}$        &4.202        &4.659    &10.89     \\
 			\\[-1em]
 			$10^{-3}$     &3.991       &4.306    &7.88   \\
 			\\[-1em]
 			$10^{-8}$      &3.246        &3.324    &2.42 \\
 			\\[-1em]
 			$10^{-12}$       & 0.339        &0.342     &1.03  \\
 			\\[-1em]
 			\bottomrule
 		\end{tabular}
 		\caption{Convergence of KL-divergence problem for different $\epsilon$.}
 		\label{tab:mc2}
\end{table}

\section{Concluding remarks}
There are two somewhat different reasons for the gap between the considered  statistical upper and deterministic lower bounds. One reason is the optimality gap similar to the risk neutral case. The additional gap, as compared to the risk neutral setting, appears because  Jensen’s inequality is employed in derivations \eqref{proofinductionub1}.  This gap tends to increase as the function $\Psi(\cdot,\theta)$ becomes more ``nonlinear". This can be clearly seen in Table \ref{tab:mc2}, the gap increases with increase of $\e$, and also in Table \ref{tab:mc} as the problem becomes more risk-averse.

When the function $\Psi$ is not polyhedral, as for instance in the setting of $\phi$-divergence example,  the procedure requires solving nonlinear optimization programs. This could be inconvenient since nonlinear optimization solvers should be used, which are known to be less efficient than linear solvers. In the considered example of KL-divergence, this requires solving one-dimensional nonlinear programs, which does not pose a significant problem. In general, in order to keep the procedure to linear programming solvers, the $Q$-factor approach, discussed in section \ref{app-q} of the Appendix, can be used. Note however that the  $Q$-factor  approach involves increasing the state space  which could significantly slow down   the convergence of the algorithm.

\bibliographystyle{plain}
\bibliography{references}

{\bf Acknowledgment}
Research of  A. Shapiro
 was partially supported by Air Force Office of Scientific Research (AFOSR)  under Grant FA9550-22-1-0244.

\setcounter{equation}{0}
\section{Appendix}
\label{append}

\subsection{Figure}

\begin{figure}[h]
\centering
\includegraphics[width=0.9\textwidth,height=0.66\textwidth]{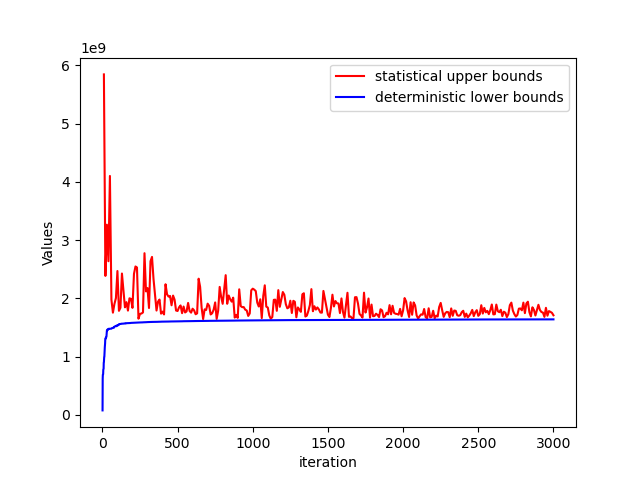}
\caption{Evolution of lower and upper bounds for convex combination of expectation and \AVaR \ problem when $\lambda = 0.5$.}
	\label{fig:evo}
\end{figure}

\subsection{Controls}
\label{app-cont}
Consider the setting where the control set depends on the state variables. That is,  consider the   extension of problem \eqref{soc} - \eqref{soc-b},
where the feasibility constraints $u_t\in \U_t$ are  replaced by $u_t\in \U_t(x_t)$ with $\U_t:\bbr^{n_t}\rr  \bbr^{m_t}$  being  a (measurable) point to set mapping, $t=1,...,T$.
By changing  the cost functions to $\bar{c}_t(x_t,u_t,\xi_t):=c_t(x_t,u_t,\xi_t)+\Ind_{\U_t(x_t)}(u_t)$, 
where  $\Ind_{\U_t(x_t)}$ is the indicator function of set $\U_t(x_t)$, 
we can write  the corresponding  problem
 in the following form
\begin{eqnarray}
\label{app-3}
&\min\limits_{\pi} &\bbe^\pi\left [ \sum_{t=1}^{T}
\bar{c}_t(x_t,u_t,\xi_t)+c_{T+1}(x_{T+1})
\right],\\
&{\rm s.t.} &
\label{app-4}
u_t=\pi_t(\xi_{[t-1]}),
u_t\in \bbr^{m_t}\; {\rm and}\;
x_{t+1}=F_t(x_t,u_t,\xi_t),\;
t=1,...,T.
\end{eqnarray}

In order to maintain convexity of the value functions,  we need to verify convexity  in $(x_t,u_t)$ of the  cost functions $\bar{c}_t(x_t,u_t,\xi_t)$, i.e., to  verify convexity of the indicator functions $\psi_t(x_t,u_t):=\Ind_{\U_t(x_t)}(u_t)$.  Note that $\psi_t(x_t,u_t)=0$ if $u_t\in \U_t(x_t)$, and $\psi_t(x_t,u_t)=+\infty$ otherwise, i.e., $\psi_t(\cdot,\cdot)$ is the indicator function of the set  $\Gr(\U_t):=\{(x_t,u_t):u_t\in \U_t(x_t)\}$     (this set is the graph of the multifunction  $\U_t$). Therefore
$\psi_t(x_t,u_t)$ is convex iff the set $\Gr(\U_t)$ is a convex subset of $\bbr^{n_t}\times\bbr^{m_t}$. In particular, suppose that
\begin{equation}\label{contset}
 \U_t(x_t):=\{u_t:g_{tk}(x_t,u_t)\le 0,\;k=1,...,K\}
\end{equation}
   for given  functions  $g_{tk}:\bbr^{n_t}\times \bbr^{m_t}\to \bbr$. Then the set
$\Gr(\U_t)$ is convex if the functions
  $g_{tk}(\cdot,\cdot)$ are  convex.

In the risk neutral  case the corresponding dynamic programming equations for the lower bounding approximations of the values functions, become
\begin{equation}\label{subd-2bb}
\underline{ V}_t(x_t)=\inf\limits_{u_t\in \U_t(x_t)}
 \sum_{i=1}^N  p_{ti} \left [c_{ti}(x_t,u_t)+
\underline{ V}_{t+1}\big(A_{ti} x_t+B_{ti} u_t+ b_{ti} \big)\right]. 
\end{equation}
Suppose that  
the set $ \U_t(x_t)$ is of the form \eqref{contset} with functions $g_{tk}(x_t,u_t)$ being convex. We need a procedure
to compute a subgradient of the right hand side of \eqref{subd-2bb}. 
Let 
\[
\underline{ V}_{t+1}(x_{t+1})=\max_{j=1,...,M}\left\{\ell_{t+1,j}
(x_{t+1})\right \}
\]
be the current representation of $\underline{ V}_{t+1}$ by its cutting planes $\ell_{t+1,j}(x_{t+1})=a_{t+1,j}^\top x_{t+1}+h_{t+1,j}$.
We can write the minimization problem \eqref{subd-2bb} as the following program
\begin{equation}\label{linprob}
\begin{array}{lll}
&\min\limits_{u,z}  &\sum_{i=1}^N  p_{ti} \left [c_{ti}^\top(x_t,u_t)+z_i\right]\\
&{\rm s.t.} & \ell_{t+1,j}(A_{ti} x_t+ B_{ti} u_t+   b_{ti} )\le z_i,\;i=1,...,N,\;
j=1,...,M,\\
&& g_{tk}(x_t,u_t)\le 0,\;k=1,...,K. 
\end{array}
\end{equation}
Suppose further  that the cost functions
$c_{ti}(x_t,u_t)$
and the constraint functions $g_{tk}(x_t,u_t)$
are   linear.  Then the above   problem \eqref{linprob} is linear. 
The required  subgradient can be computed by solving the dual of the linear program \eqref {linprob}.

In the risk averse case it is possible to proceed in a similar way.  Suppose for example $\R_t=\avr_\alpha$ risk measure. Then we can write the corresponding dynamic equations in the form
 \begin{equation}\label{subd-2cc}
\underline{ V}_t(x_t)=\inf\limits_{u_t\in \U_t(x_t),\,\theta\in \bbr}\left\{ \theta+\alpha^{-1}    \sum_{i=1}^N p_{ti} 
\left[    c_{ti}(x_t,u_t)+
\underline{ V}_{t+1}\big(A_{ti} x_t+B_{ti} u_t+ b_{ti} \big) -\theta \right]_+ \right\}.
\end{equation}
In the above formulation controls and  parameter $\theta$  of the $\avr_\alpha$ risk measure are computed simultaneously.   The minimization problem \eqref{subd-2cc} can be written as the following program
\begin{eqnarray*}
&\min\limits_{u,\theta,z}& \theta+\alpha^{-1}  \sum_{i=1}^N   p_{ti} z_i\\
&{\rm s.t.} & c_{ti}(x_t,u_t)+
\ell_{t+1,j}(A_{ti} x_t+ B_{ti} u_t+   b_{ti} )-\theta\le z_i,\;i=1,...,N,\;
j=1,...,M,  \\
&&0\le z_i,\;i=1,...,N,\\
&& g_{tk}(x_t,u_t)\le 0,\;k=1,...,K. 
\end{eqnarray*}
If the cost functions
$c_{ti}(x_t,u_t)$
and the constraint functions $g_{tk}(x_t,u_t)$ are linear, this is a linear program. In general it is possible to write problem \eqref{subd-2cc} as a linear program if the risk measure and  the cost functions are polyhedral and the constraint functions are linear.

\subsection{Optimal Control and Stochastic Programming modeling }
\label{app-stoch}
Mainly for historical reasons,  the SDDP algorithm was formulated first in the framework of the SP  modeling. Quite often the same  optimization problem can be alternatively formulated either in the SOC or SP framework.
In both cases the decision should be based on information available at time of the decision, this is the so-called nonaticipativity principle. There are various ways how the information available at time $t$ can be represented. Here we assume that it is defined by history of the random (data) process $\xi_t$. We label the available   history at time $t$ as $\xi_{[t-1]}=(\xi_0,\xi_1,...,\xi_{t-1})$, with $\xi_0$ being given (deterministic). Of course,  shifting the time   label we can write  this as $\xi_{[t]}=(\xi_1,...,\xi_{t})$  with now  $\xi_1$ being deterministic representing the initial conditions,   which is more common in the SP framework. What is important is that in both cases our decisions are functions of the observed realizations of the data process at time of the decision. It also could be noted that we need to  consider only policies which are  functions of the data process alone because of the basic assumption that the distribution of the random process $\xi_t$ does not depend on our decisions.

One important difference between the SOC and SP modeling is that in the SOC approach there is a clear separation between the states and controls. Because of the stagewise independence assumption, the value functions $V_t(x_t)$ are functions of the state variables only. The controls $u_t$ and the corresponding values $\theta_t$ of the parameter vector are computed (estimated) simultaneously based on equation  \eqref{risk-00a}. That is, the estimated values of $\theta_t$ are functions of state $x_t$ and optimal controls $\bar{u}_t$,   based on a current approximation of the value function (see eq. \eqref{minimizer}).
This makes the computed estimates of $\theta_t$  to be consistent for the generated discretization  (sample) of the marginal distribution of $\xi_t$. This is in contrast to the SP approach where the bias of the corresponding estimates of $\theta_t$ explodes  exponentially with increase of the number of stages (cf., \cite{ShapiroDing}).

\subsection{$Q$-factor approach }
\label{app-q}

The following is a counterpart of  the $Q$-factor approach popular in the SOC applications.
Consider the dynamic equations 
 \begin{equation}
 \label{risk-aa}
 V_t(x_t) =
 \inf\limits_{u_t\in \U_t,\,\theta_t\in \Theta}  \bbe_{P_t} \left[ \Psi\big (c_{t}(x_t,u_t,\xi_t)+ V_{t+1} (A_{t} x_t+B_{t} u_t+ b_{t}),\theta_t  \big )\right],
 \end{equation}
  and define
\begin{equation}\label{qfact-1}
Q_t(x_t,u_t,\theta_t):=\bbe_{P_t} \left[ \Psi\big (c_{t}(x_t,u_t,\xi_t)+ V_{t+1} (A_{t} x_t+B_{t} u_t+ b_{t}),\theta_t  \big )\right].
\end{equation}
We have that
\[
 V_t(x_t)=
 \inf\limits_{u_t\in \U_t,\,\theta_t\in \Theta}   Q_t(x_t,u_t,\theta_t),
 \]
 and hence the dynamic equations \eqref{risk-aa} can be written in terms of $Q_t(x_t,u_t,\theta_t)$ as
\begin{equation}\label{qfact-2}
Q_t(x_t,u_t,\theta_t)=\bbe_{P_t} \Big[ \Psi  \Big(c_{t}(x_t,u_t,\xi_t)+
 \inf\limits_{u_{t+1}\in \U_{t+1},\,\theta_{t+1}\in \Theta}   Q_{t+1}\big( A_{t} x_t+B_{t} u_t+ b_{t},   u_{t+1},\theta_{t+1} \big),\theta_t\Big)\Big].
\end{equation}

The cutting planes,  SDDP type,  algorithm can be applied directly to  functions $Q_t(x_t,u_t,\theta_t)$ rather than to the value functions $V_t(x_t)$. In the backward step of the algorithm,  subgradients with respect to $x_t,u_t$ and $\theta_t$,   of the current approximations of the functions $Q_t(x_t,u_t,\theta_t)$,  should be computed.
An advantage of that approach is that the calculation of these subgradients does not require
solving {\em nonlinear optimization} programs  even if the function $\Psi$ is not polyhedral\footnote{The    function $\Psi$ is not polyhedral, for example,
in the $\phi$-divergence case. In that case the SDDP algorithm, applied to the value functions $V_t(x_t)$,   requires solving nonlinear programs.}.
On the other hand,  this $Q$-factor  approach involves increasing the state space from $x_t$  to $(x_t,u_t,\theta_t)$, which could make the convergence of the algorithm considerably slower.

\end{document}